\documentclass[a4paper,twoside,11pt]{article}
\usepackage{amsfonts}
\usepackage{amssymb}\usepackage{amsmath}
\newtheorem{thm}{Theorem}[section]
 
 \numberwithin{equation}{section}
\newcommand{\double}{\baselineskip 1.24 \baselineskip}

\title{Fourier Cosine and Sine Transform on fractal space
}

\author{{Guang-Sheng  Chen\thanks{\text{E-mail address}: cgswavelets@126.com(Chen)
}\quad}\\
{\small Department of Computer Engineering, Guangxi Modern
Vocational Technology College,} \\{\small Hechi,Guangxi, 547000,
P.R. China}
}
\begin{document}
\date{}
\maketitle \double

\textbf{Abstract:}\quad In this paper, we establish local fractional Fourier
Cosine and Sine Transforms on fractal space, considered some properties of
local fractional Cosine and Sine Transforms, show applications of local
fractional Fourier Cosine and Sine transform to local fractional equations
with local fractional derivative.  

\section{Introduction}
\hskip\parindent%缩进
The fractional Fourier transform has been investigated in a number of papers
and has been proved to be very useful in solving engineering problems [1-7].
It is important to deal with the continuous fractal functions, which are
irregular in the real world. Recently, Yang-Fourier transform based on the
local fractional calculus was introduced [8] and Yang continued to study
this subject [9-10]. The importance of Yang-Fourier transform for fractal
functions derives from the fact that this is the only mathematic model which
focuses on local fractional continuous functions derived from local
fractional calculus. More recently, some model for engineering derived from
local fractional derivative was proposed [11-12].

The purpose of this paper is to establish local fractional Cosine and Sine
Transforms based on the Yang-Fourier transforms and consider its application
to local fractional equations with local fractional derivative. This paper
is organized as follows. In section 2, local fractional Cosine and Sine
Transforms is derived; Section 3 presents Properties of local fractional
Fourier Cosine and Sine Transforms; Applications of local fractional Fourier
Cosine and Sine Transforms are discussed in section 4.

\section{local fractional Cosine and Sine Transforms}
\hskip\parindent%缩进
In this section, we start with the following result [9,10]:

\begin{equation*}
f(x) = \frac{1}{(2\pi )^\alpha }\int_{ - \infty }^\infty {C_k E_\alpha
(i^\alpha x^\alpha \omega ^\alpha )(d\omega )^\alpha },
\tag{2.1}
\label{2.1}
\end{equation*}
where
\begin{equation*}
C_k = \frac{1}{\Gamma (1 + \alpha )}\int_{ - \infty }^\infty
{f(x)E_\alpha ( - i^\alpha x^\alpha \omega ^\alpha )(dx)^\alpha }.
\tag{2.2}
\label{2.2}
\end{equation*}
From (2.2) , the Yang-Fourier transform of $f(x)$ is given by [9,10]
\begin{equation*}
F_\alpha \{f(x)\} = f_\omega ^{F,\alpha } (\omega ): = \frac{1}{\Gamma (1 +
\alpha )}\int_{ - \infty }^\infty {E_\alpha ( - i^\alpha \omega ^\alpha
x^\alpha )f(x)(dx)^\alpha }.
\tag{2.3}
\label{2.3}
\end{equation*}
And its Inverse formula of Yang-Fourier's transforms as follows
\begin{equation*}
f(x) = F_\alpha ^{ - 1} (f_\omega ^{F,\alpha } (\omega )): = \frac{1}{(2\pi
)^\alpha }\int_{ - \infty }^\infty {E_\alpha (i^\alpha \omega ^\alpha
x^\alpha )f_\omega ^{F,\alpha } (\omega )(d\omega )^\alpha }.
\tag{2.4}
\label{2.4}
\end{equation*}
Now ,by (2.1) and (2.2) ,we have
\begin{equation*}
f(x) = \frac{1}{(2\pi )^\alpha }\int_{ - \infty }^\infty {[\frac{1}{\Gamma
(1 + \alpha )}\int_{ - \infty }^\infty {f(\xi )E_\alpha ( - i^\alpha \xi
^\alpha \omega ^\alpha )(d\xi )^\alpha } ]E_\alpha (i^\alpha x^\alpha \omega
^\alpha )(d\omega )^\alpha }.
\tag{2.5}
\label{2.5}
\end{equation*}
Here, we named (2.5) the Yang-Fourier integral formula.

We express the exponential factor $E_\alpha (i^\alpha \omega ^\alpha
(x^\alpha - \xi ^\alpha ))$ in (2.5) in terms of trigonometric functions on
fractal set of fractal dimension$\alpha $and use the even and odd nature of
the cosine and the sine functions respectively as functions of $\omega $, so
that (2.5) can be written as
\begin{equation*}
\begin{split}
& f(x) = \frac{1}{(2\pi )^\alpha }\int_{ - \infty }^\infty {[\frac{1}{\Gamma
(1 + \alpha )}\int_{ - \infty }^\infty {f(\xi )E_\alpha ( - i^\alpha \xi
^\alpha \omega ^\alpha )(d\xi )^\alpha } ]E_\alpha (i^\alpha x^\alpha \omega
^\alpha )(d\omega )^\alpha } \\
 &= \frac{1}{(2\pi )^\alpha }\int_{ - \infty }^\infty {[\frac{1}{\Gamma (1 +
\alpha )}\int_{ - \infty }^\infty {f(\xi )E_\alpha (i^\alpha \omega ^\alpha
(x^\alpha - \xi ^\alpha ))(d\xi )^\alpha } ](d\omega )^\alpha } \\
 &= \frac{2}{(2\pi )^\alpha }\int_0^\infty {[\frac{1}{\Gamma (1 + \alpha
)}\int_{ - \infty }^\infty {f(\xi )\cos _\alpha (\omega ^\alpha (x^\alpha -
\xi ^\alpha ))(d\xi )^\alpha } ](d\omega )^\alpha }. \\
\end{split}
\notag
\end{equation*}
Hence,we have
\begin{equation*}
f(x) = \frac{2}{(2\pi )^\alpha }\int_0^\infty {[\frac{1}{\Gamma (1 + \alpha
)}\int_{ - \infty }^\infty {f(\xi )\cos _\alpha (\omega ^\alpha (x^\alpha -
\xi ^\alpha ))(d\xi )^\alpha } ](d\omega )^\alpha }.
\tag{2.6}
\label{2.6}
\end{equation*}
This is another version of the Yang-Fourier integral formula.

We now assume that $f(x)$ is an even function and expand the cosine function
in (2.6) to obtain
\begin{equation*}
f(x) = f( - x) = \frac{4}{(2\pi )^\alpha }\frac{1}{\Gamma (1 + \alpha
)}\int_0^\infty {\cos _\alpha (\omega ^\alpha x^\alpha )(d\omega )^\alpha }
\int_0^\infty {f(\xi )\cos _\alpha (\omega ^\alpha \xi ^\alpha )(d\xi
)^\alpha }.
\tag{2.7}
\label{2.7}
\end{equation*}
This is called the local fractional Fourier cosine integral formula.

Similarly, for an odd function $f(x)$, we obtain the local fractional
Fourier sine integral formula
\begin{equation*}
f(x) = - f( - x) = \frac{4}{(2\pi )^\alpha }\frac{1}{\Gamma (1 + \alpha
)}\int_0^\infty {\sin _\alpha (\omega ^\alpha x^\alpha )(d\omega )^\alpha }
\int_0^\infty {f(\xi )\sin _\alpha (\omega ^\alpha \xi ^\alpha )(d\xi
)^\alpha }. \tag{2.8}
\label{2.8}
\end{equation*}
The local fractional Fourier cosine integral formula (2.7) leads to the
local fractional Fourier cosine transform and its inverse defined by
\begin{equation*}
F_{\alpha ,c} \{f(x)\} = f_{\omega ,c}^{F,\alpha } (\omega ): =
\frac{2}{\Gamma (1 + \alpha )}\int_0^\infty {f(x)\cos _\alpha (\omega
^\alpha x^\alpha )(dx)^\alpha },  \tag{2.9}
\label{2.9}
\end{equation*}
\begin{equation*}
f(x) = F_{\alpha .c}^{ - 1} (f_{\omega ,c}^{F,\alpha } (\omega )): =
\frac{2}{(2\pi )^\alpha }\int_0^\infty {f_{\omega ,c}^{F,\alpha } (\omega
)\cos _\alpha (\omega ^\alpha x^\alpha )(d\omega )^\alpha }, \tag{2.10}
\label{2.10}
\end{equation*}
where $F_{\alpha ,c} $ is the local fractional Fourier cosine transform
operator and $F_{\alpha .c}^{ - 1} $ is its inverse operator.

Similarly, the local fractional Fourier sine integral formula (2.8) leads to
the local fractional Fourier sine transform and its inverse defined by
\begin{equation*}
F_{\alpha ,s} \{f(x)\} = f_{\omega ,s}^{F,\alpha } (\omega ): =
\frac{2}{\Gamma (1 + \alpha )}\int_0^\infty {f(x)\sin _\alpha (\omega
^\alpha x^\alpha )(dx)^\alpha },  \tag{2.11}
\label{2.11}
\end{equation*}
\begin{equation*}
f(x) = F_{\alpha .s}^{ - 1} (f_{\omega ,s}^{F,\alpha } (\omega )): =
\frac{2}{(2\pi )^\alpha }\int_0^\infty {f_{\omega ,c}^{F,\alpha } (\omega
)\sin _\alpha (\omega ^\alpha x^\alpha )(d\omega )^\alpha },\tag{2.12}
\label{2.12}
\end{equation*}
where $F_{\alpha ,s} $ is the local fractional Fourier cosine transform
operator and $F_{\alpha .s}^{ - 1} $ is its inverse operator.
\\ {\bf Example 1.} Show that $a > 0$
\begin{equation*}
F_{\alpha ,s} \{E_\alpha ( - at^\alpha )\} = \frac{2\omega ^\alpha }{a^2 +
\omega ^{2\alpha }},
\quad
a > 0, \tag{2.13}
\label{2.13}
\end{equation*}
\begin{equation*}
F_{\alpha ,c} \{E_\alpha ( - at^\alpha )\} = \frac{2a}{a^2 + \omega
^{2\alpha }}, \quad
a > 0. \tag{2.14}
\label{2.14}
\end{equation*}
By local fractional Fourier cosine transform,we have
\begin{equation*}
\begin{split}
 &F_{\alpha ,s} \{E_\alpha ( - at^\alpha )\} = \frac{2}{\Gamma (1 + \alpha
)}\int_0^\infty {E_\alpha ( - at^\alpha )\sin _\alpha (\omega ^\alpha
x^\alpha )(dx)^\alpha } \\
 &= \frac{2\omega ^\alpha }{a^2} - \frac{2\omega ^{2\alpha
}}{a^2}\frac{1}{\Gamma (1 + \alpha )}\int_0^\infty {E_\alpha ( - at^\alpha
)\sin _\alpha (\omega ^\alpha x^\alpha )(dx)^\alpha } \\
 &= \frac{2\omega ^\alpha }{a^2} - \frac{\omega ^{2\alpha }}{a^2}F_{\alpha
,s} \{E_\alpha ( - at^\alpha )\}. \\
\end{split}
\tag{2.15}
\label{2.15}
\end{equation*}
By (2.15),we obtain
\[
F_{\alpha ,s} \{E_\alpha ( - at^\alpha )\} = \frac{2\omega ^\alpha }{a^2 +
\omega ^{2\alpha }}.
\]
Similarly, we obtain (2.14)
\section{Properties of local fractional Fourier Cosine and Sine Transforms}
\begin{thm} \label{THEOREM 1} If $F_{\alpha ,c} \{f(x)\} = f_{\omega ,c}^{F,\alpha }
(\omega )$ and $F_{\alpha ,s} \{f(x)\} = f_{\omega ,s}^{F,\alpha } (\omega
)$, then
\begin{equation*}
F_{\alpha ,c} \{f(ax)\} = \frac{1}{a^\alpha }f_{\omega ,c}^{F,\alpha }
\left( {\frac{\omega }{a}} \right),
\tag{3.1}
\label{3.1}
\end{equation*}
\begin{equation*}
F_{\alpha ,s} \{f(ax)\} = \frac{1}{a^\alpha }f_{\omega ,s}^{F,\alpha }
\left( {\frac{\omega }{a}} \right).
\tag{3.2}
\label{3.2}
\end{equation*}
\end{thm}
\textbf{Proof}. As a direct application of the local fractional Fourier
Cosine transform, we derive the follow identity
\begin{equation*}
\begin{split}
& F_{\alpha ,c} \{f(ax)\} = \frac{2}{\Gamma (1 + \alpha )}\int_0^\infty
{f(ax)\cos _\alpha (\omega ^\alpha x^\alpha )(dx)^\alpha } \\
& = \frac{2}{a^\alpha \Gamma (1 + \alpha )}\int_0^\infty {f(ax)\cos _\alpha
(\omega ^\alpha x^\alpha )(dax)^\alpha }. \\
 \end{split}
\tag{3.3}
\label{3.3}
\end{equation*}
\noindent
taking $y = ax$ in (3.3) implies that
\begin{equation*}
\begin{split}
 &\frac{2}{a^\alpha \Gamma (1 + \alpha )}\int_0^\infty {f(ax)\cos _\alpha
(\omega ^\alpha x^\alpha )(dax)^\alpha } \\
 &= \frac{2}{a^\alpha \Gamma (1 + \alpha )}\int_0^\infty {f(ax)\cos _\alpha
((\textstyle{\omega \over a})^\alpha (ax)^\alpha )(dax)^\alpha } \\
& = \frac{2}{a^\alpha \Gamma (1 + \alpha )}\int_0^\infty {f(y)\cos _\alpha
((\textstyle{\omega \over a})^\alpha y^\alpha )(dy)^\alpha } =
\frac{1}{a^\alpha }f_{\omega ,c}^{F,\alpha } \left( {\frac{\omega }{a}}
\right). \\
 \end{split}
\notag
\end{equation*}
Similarly, we obtain (3.2)

Under appropriate conditions, the following properties also hold:
\begin{equation*}
F_{\alpha ,c} \{f^{(\alpha )}(x)\} = \omega ^\alpha f_{\omega ,c}^{F,\alpha
} (\omega ) - 2f(0),
\tag{3.4}
\label{3.4}
\end{equation*}
\begin{equation*}
F_{\alpha ,c} \{f^{(2\alpha )}(x)\} = - \omega ^{2\alpha }f_{\omega
,c}^{F,\alpha } (\omega ) - 2f^{(\alpha )}(0),
\tag{3.5}
\label{3.5}
\end{equation*}
\begin{equation*}
F_{\alpha ,s} \{f^{(\alpha )}(x)\} = - \omega ^\alpha f_{\omega
,c}^{F,\alpha } (\omega ),
\tag{3.6}
\label{3.6}
\end{equation*}
\begin{equation*}
F_{\alpha ,s} \{f^{(2\alpha )}(x)\} = - \omega ^{2\alpha }f_{\omega
,c}^{F,\alpha } (\omega ) + 2\omega ^\alpha f(0),
\tag{3.7}
\label{3.7}
\end{equation*}
These results can be generalized for the cosine and sine transforms of
higher order derivatives of a function.
\begin{thm}(Convolution Theorem for the local fractional Fourier
Cosine Transform)\label{THEOREM 2(Convolution Theorem for the local fractional Fourier
Cosine Transform)} If $F_{\alpha ,c} \{f(x)\} = f_{\omega ,c}^{F,\alpha } (\omega )$ and
$F_{\alpha ,c} \{g(x)\} = g_{\omega ,c}^{F,\alpha } (\omega )$,  then
\begin{equation*}
\begin{split}
& \frac{2}{(2\pi )^\alpha }\int_0^\infty {f_{\omega ,c}^{F,\alpha } (\omega
)g_{\omega ,c}^{F,\alpha } (\omega )\cos _\alpha (\omega ^\alpha x^\alpha
)(d\omega )^\alpha } \\
 &= \frac{1}{\Gamma (1 + \alpha )}\int_0^\infty {f(\xi )[g(x + \xi ) +
g(\vert x - \xi \vert )](d\xi )^\alpha }. \\
\end{split}
\tag{3.8}
\label{3.8}
\end{equation*}
Or, equivalently,
\begin{equation*}
\begin{split}
& \int_0^\infty {f_{\omega ,c}^{F,\alpha } (\omega )g_{\omega ,c}^{F,\alpha }
(\omega )\cos _\alpha (\omega ^\alpha x^\alpha )(d\omega )^\alpha } \\
 &= \frac{(2\pi )^\alpha }{2\Gamma (1 + \alpha )}\int_0^\infty {f(\xi )[g(x +
\xi ) + g(\vert x - \xi \vert )](d\xi )^\alpha }.\\
 \end{split}
\tag{3.9}
\label{3.9}
\end{equation*}
\end{thm}
\textbf{Proof.}Using the definition of the inverse local fractional Fourier
cosine transform, we Have
\begin{equation*}
\begin{split}
 &F_{\alpha .c}^{ - 1} (f_{\omega ,c}^{F,\alpha } (\omega )g_{\omega
,c}^{F,\alpha } (\omega )) = \frac{2}{(2\pi )^\alpha }\int_0^\infty
{f_{\omega ,c}^{F,\alpha } (\omega )g_{\omega ,c}^{F,\alpha } (\omega )\cos
_\alpha (\omega ^\alpha x^\alpha )(d\omega )^\alpha } \\
 &= \frac{2}{(2\pi )^\alpha }\int_0^\infty {g_{\omega ,c}^{F,\alpha } (\omega
)\cos _\alpha (\omega ^\alpha x^\alpha )(d\omega )^\alpha } \frac{2}{\Gamma
(1 + \alpha )}\int_0^\infty {f(x)\cos _\alpha (\omega ^\alpha x^\alpha
)(dx)^\alpha } \\
 &= \frac{4}{(2\pi )^\alpha \Gamma (1 + \alpha )}\int_0^\infty {g_{\omega
,c}^{F,\alpha } (\omega )\cos _\alpha (\omega ^\alpha x^\alpha )(d\omega
)^\alpha } \int_0^\infty {f(x)\cos _\alpha (\omega ^\alpha x^\alpha
)(dx)^\alpha }. \\
\end{split}
\notag
\end{equation*}
Hence,we obtain
\begin{equation*}
\begin{split}
 &F_{\alpha .c}^{ - 1} (f_{\omega ,c}^{F,\alpha } (\omega )g_{\omega
,c}^{F,\alpha } (\omega )) \\
 &= \frac{4}{(2\pi )^\alpha \Gamma (1 + \alpha )}\int_0^\infty {g_{\omega
,c}^{F,\alpha } (\omega )\cos _\alpha (\omega ^\alpha x^\alpha )(d\omega
)^\alpha } \int_0^\infty {f(\xi )\cos _\alpha (\omega ^\alpha \xi ^\alpha
)(d\xi )^\alpha } \\
& = \frac{4}{(2\pi )^\alpha \Gamma (1 + \alpha )}\int_0^\infty {f(\xi )(d\xi
)^\alpha } \int_0^\infty {g_{\omega ,c}^{F,\alpha } (\omega )\cos _\alpha
(\omega ^\alpha x^\alpha )\cos _\alpha (\omega ^\alpha \xi ^\alpha )(d\omega
)^\alpha } \\
 &= \frac{2}{(2\pi )^\alpha \Gamma (1 + \alpha )}\int_0^\infty {f(\xi )(d\xi
)^\alpha } \int_0^\infty {g_{\omega ,c}^{F,\alpha } (\omega )[\cos _\alpha
\omega ^\alpha (x^\alpha + \xi ^\alpha ) + \cos _\alpha \omega ^\alpha
(\vert x^\alpha - \xi ^\alpha \vert )](d\omega )^\alpha } \\
 &= \frac{1}{\Gamma (1 + \alpha )}\int_0^\infty {f(\xi )[g(x + \xi ) +
g(\vert x - \xi \vert )](d\xi )^\alpha }. \\
\end{split}
\notag
\end{equation*}
in which the definition of the inverse local fractional Fourier cosine
transform is used. This proves (3.8). It also follows from the proof of
Theorem 2 that
\[
\int_0^\infty {f_{\omega ,c}^{F,\alpha } (\omega )g_{\omega ,c}^{F,\alpha }
(\omega )\cos _\alpha (\omega ^\alpha x^\alpha )(d\omega )^\alpha } =
\frac{(2\pi )^\alpha }{2\Gamma (1 + \alpha )}\int_0^\infty {f(\xi )[g(x +
\xi ) + g(\vert x - \xi \vert )](d\xi )^\alpha }.
\]
This proves result (3.9).

Putting $x = 0$ in (3.9), we obtain
\[
\int_0^\infty {f_{\omega ,c}^{F,\alpha } (\omega )g_{\omega ,c}^{F,\alpha }
(\omega )(d\omega )^\alpha } = \frac{(2\pi )^\alpha }{\Gamma (1 + \alpha
)}\int_0^\infty {f(\xi )g(\xi )(d\xi )^\alpha } = \frac{(2\pi )^\alpha
}{\Gamma (1 + \alpha )}\int_0^\infty {f(x)g(x)(dx)^\alpha }.
\]
Substituting $g(x) = \overline {f(x)} $ gives, since $f_{\omega
,c}^{F,\alpha } (\omega ) = \overline {g_{\omega ,c}^{F,\alpha } (\omega )}
$,
\begin{equation*}
\int_0^\infty {\vert f_{\omega ,c}^{F,\alpha } (\omega )\vert ^2(d\omega
)^\alpha } = \frac{(2\pi )^\alpha }{\Gamma (1 + \alpha )}\int_0^\infty
{\vert f(x)\vert ^2(dx)^\alpha } .
\tag{3.10}
\label{3.10}
\end{equation*}
This is the Parseval relation for the local fractional Fourier cosine
transform.

Similarly, we obtain
\begin{equation*}
\begin{split}
 &\int_0^\infty {f_{\omega ,s}^{F,\alpha } (\omega )g_{\omega ,s}^{F,\alpha }
(\omega )\cos _\alpha (\omega ^\alpha x^\alpha )(d\omega )^\alpha } \\
 &= \int_0^\infty {f_{\omega ,s}^{F,\alpha } (\omega )g_{\omega ,s}^{F,\alpha
} (\omega )\cos _\alpha (\omega ^\alpha x^\alpha )(d\omega )^\alpha } \\
 &= \frac{2}{\Gamma (1 + \alpha )}\int_0^\infty {g_{\omega ,s}^{F,\alpha }
(\omega )\cos _\alpha (\omega ^\alpha x^\alpha )(d\omega )^\alpha }
\int_0^\infty {f(\xi )\sin _\alpha (\omega ^\alpha \xi ^\alpha )(d\xi
)^\alpha }. \\
 \end{split}
\notag
\end{equation*}
which is, by interchanging the order of integration,
\begin{equation*}
\begin{split}
 &\int_0^\infty {f_{\omega ,s}^{F,\alpha } (\omega )g_{\omega ,s}^{F,\alpha }
(\omega )\cos _\alpha (\omega ^\alpha x^\alpha )(d\omega )^\alpha } \\
 &= \frac{2}{\Gamma (1 + \alpha )}\int_0^\infty {g_{\omega ,s}^{F,\alpha }
(\omega )\cos _\alpha (\omega ^\alpha x^\alpha )(d\omega )^\alpha }
\int_0^\infty {f(\xi )\sin _\alpha (\omega ^\alpha \xi ^\alpha )(d\xi
)^\alpha } \\
 &= \frac{2}{\Gamma (1 + \alpha )}\int_0^\infty {f(\xi )(d\xi )^\alpha }
\int_0^\infty {g_{\omega ,s}^{F,\alpha } (\omega )\sin _\alpha (\omega
^\alpha \xi ^\alpha )\cos _\alpha (\omega ^\alpha x^\alpha )(d\omega
)^\alpha } \\
 &= \frac{1}{\Gamma (1 + \alpha )}\int_0^\infty {f(\xi )(d\xi )^\alpha }
\int_0^\infty {g_{\omega ,s}^{F,\alpha } (\omega )[\sin _\alpha \omega
^\alpha (x^\alpha + \xi ^\alpha ) + \sin _\alpha \omega ^\alpha (\xi ^\alpha
- x^\alpha )](d\omega )^\alpha } \\
 &= \frac{(2\pi )^\alpha }{2\Gamma (1 + \alpha )}\int_0^\infty {f(\xi )[g(\xi
+ x) + g(\xi - x)](d\xi )^\alpha }. \\
\end{split}
\notag
\end{equation*}
in which the inverse local fractional Fourier sine transform is used. Thus,
we find
\begin{equation*}
\begin{split}
&\int_0^\infty {f_{\omega ,s}^{F,\alpha } (\omega )g_{\omega ,s}^{F,\alpha }
(\omega )\cos _\alpha (\omega ^\alpha x^\alpha )(d\omega )^\alpha }\\
& =
\frac{(2\pi )^\alpha }{2\Gamma (1 + \alpha )}\int_0^\infty {f(\xi )[g(\xi +
x) + g(\xi - x)](d\xi )^\alpha }.
\end{split}
\tag{3.11}
\label{3.11}
\end{equation*}
Or, equivalently,
\begin{equation*}
\begin{split}
 &\frac{2}{(2\pi )^\alpha }\int_0^\infty {f_{\omega ,s}^{F,\alpha } (\omega
)g_{\omega ,s}^{F,\alpha } (\omega )\cos _\alpha (\omega ^\alpha x^\alpha
)(d\omega )^\alpha } \\
 &= \frac{1}{\Gamma (1 + \alpha )}\int_0^\infty {f(\xi )[g(\xi + x) + g(\xi -
x)](d\xi )^\alpha }.\\
 \end{split}
\tag{3.12}
\label{3.12}
\end{equation*}
Result (3.11) or (3.12) is also called the Convolution Theorem of the local
fractional Fourier cosine transform.

Putting $x = 0$ in (3.11) gives
\[
\int_0^\infty {f_{\omega ,s}^{F,\alpha } (\omega )g_{\omega ,s}^{F,\alpha }
(\omega )(d\omega )^\alpha } = \frac{(2\pi )^\alpha }{\Gamma (1 + \alpha
)}\int_0^\infty {f(\xi )g(\xi )(d\xi )^\alpha } = \frac{(2\pi )^\alpha
}{\Gamma (1 + \alpha )}\int_0^\infty {f(x)g(x)(dx)^\alpha }.
\]
Replacing $g(x))$ by $\overline {f(x)} $ gives the Parseval relation for the
local fractional Fourier sine transform
\begin{equation*}
\int_0^\infty {\vert f_{\omega ,s}^{F,\alpha } (\omega )\vert ^2(d\omega
)^\alpha } = \frac{(2\pi )^\alpha }{\Gamma (1 + \alpha )}\int_0^\infty
{\vert f(x)\vert ^2(dx)^\alpha }.
\tag{3.13}
\label{3.13}
\end{equation*}

\section{Applications}
\hskip\parindent%缩进 
Use the local fractional Fourier sine transform to solve the following
differential equation:
\[
y^{(2\alpha )} - 9y(t) = 50E_\alpha ( - 2t^\alpha ),
\quad
y(0) = y_0.
\]

Since we are interested in positive + region, we can take $y(t)$ to be an
odd function and take local fractional Fourier sine transforms. It is clear
from its definition that local fractional Fourier sine transform is linear
\[
F_{\alpha ,s} \{af_1 (x) + bf_2 (x)\} = aF_{\alpha ,s} \{f_1 (x)\} +
bF_{\alpha ,s} \{f_2 (x)\}.
\]
Using this property and taking local fractional Fourier sine transform of
both sides of the differential equation, we have
\[
F_{\alpha ,s} \{y^{(2\alpha )}\} - 9F_{\alpha ,s} \{y(t)\} = 50F_{\alpha ,s}
\{E_\alpha ( - 2t^\alpha )\}.
\]
Since
\[
F_{\alpha ,s} \{y^{(2\alpha )}(t)\} = - \omega ^{2\alpha }F_{\alpha ,s}
\{y(t)\} + 2\omega ^\alpha y(0).
\]
So
\[
 - \omega ^{2\alpha }F_{\alpha ,s} \{y(t)\} + 2\omega ^\alpha y(0) -
9F_{\alpha ,s} \{y(t)\} = 50F_{\alpha ,s} \{E_\alpha ( - 2t^\alpha )\}.
\]
which, after collecting terms, becomes
\[
(\omega ^{2\alpha } + 9)F_{\alpha ,s} \{y(t)\} = - 50\frac{2\omega ^\alpha
}{4 + \omega ^{2\alpha }} + 2\omega ^\alpha y_0.
\]
Thus
\[
F_{\alpha ,s} \{y(t)\} = - 50\frac{2\omega ^\alpha }{4 + \omega ^{2\alpha
}}\frac{1}{\omega ^{2\alpha } + 9} + y_0 \frac{2\omega ^\alpha }{\omega
^{2\alpha } + 9}.
\]
With partial fraction of
\[
\frac{1}{(4 + \omega ^{2\alpha })(\omega ^{2\alpha } + 9)} = \frac{1}{5(4 +
\omega ^{2\alpha })} - \frac{1}{5(\omega ^{2\alpha } + 9)}.
\]
we have
\begin{equation*}
\begin{split}
& F_{\alpha ,s} \{y(t)\} = - 10\frac{2\omega ^\alpha }{(4 + \omega ^{2\alpha
})} + 10\frac{2\omega ^\alpha }{5(\omega ^{2\alpha } + 9)}] + y_0
\frac{2\omega ^\alpha }{\omega ^{2\alpha } + 9} \\
 &= (y_0 + 10)\frac{2\omega ^\alpha }{\omega ^{2\alpha } + 9} -
10\frac{2\omega ^\alpha }{(4 + \omega ^{2\alpha })} \\
 &= (y_0 + 10)F_{\alpha ,s} \{E_\alpha ( - 3t^\alpha )\} - 10F_{\alpha ,s}
\{E_\alpha ( - 2t^\alpha )\}. \\
 \end{split}
 \notag
\end{equation*}
Taking the inverse local fractional Fourier sine transform, we get the
solution
\[
y(t) = (y_0 + 10)E_\alpha ( - 3t^\alpha ) - 10E_\alpha ( - 2t^\alpha ).
\]

 \end{document}